\documentclass[11pt]{amsart}

\usepackage{import}

\usepackage{set}

\usepackage[T1]{fontenc}

 \date{\today}
 

\title{Profinite properties of residually free groups}
\author{Ismael Morales}

\begin{document}
\maketitle

\begin{abstract} Wilton classified when a prime three-manifold $M$ has a residually free fundamental group $\pi_1 M$. We prove that the groups $\pi_1 M\times \Z^n$ are profinitely rigid within finitely generated residually free groups. We also establish other profinite invariants of the class of residually free groups such as coherence and subgroup separability. In the course of our proofs, we generalise a lemma of Wilton and Zalesskii on profinitely recognising when a central extension of groups splits. 
\end{abstract}

\section{Introduction}
The  class of residually free groups has attracted much attention in geometric group theory after its connection with logic explored by Sela \cite{Sel01}. It has also been a source of inspiration in questions about profinite rigidity, which become more tractable when restricting to this class. An elementary invariant of a group $G$ one usually studies to tackle these sorts of questions is its {\it $p$-cohomological dimension}, denoted by $\cd_p(G)$. This is defined as the smallest integer $m$ such that $H^k(H; \p)=0$ for all $k>m$ and for all finite-index subgroups $H<G$. If $G$ is good in the sense of Serre then $\cd_p(G)=\cd_p(\hat{G})$, hence $\cd_p$ is a profinite invariant among good groups. This invariant has been extensively studied to recognise important classes of groups. For instance, Wilton showed in  \cite{Wil18} that a limit group $L$ has $ \cd_p(L)=1$ if and only if $L$ is free. This confirmed Remeslennikov conjecture on the profinite rigidity of free groups for the class of limit groups (which are good by the work of  Grunewald--Jaikin-Zapirain--Zalesskii \cite{Gru08}). 

Moving on to dimension two, the previous was followed by \cite{Wil21}, where it is shown that a limit group $L$ has $\PD_2$ profinite completion at some prime $p$ (in the sense of \cite[Definition 3.7.1]{Neu08})   only if $L$ is a surface group itself. By a {\it surface group} we will always mean the fundamental group of a closed connected surface. The latter result was proven with an alternative argument by Fruchter and the author \cite[Corollary E]{FruchterMorales24}, where we discuss a much broader picture about the recognition of free and surface groups from Betti numbers of finite-index subgroups. With entirely different methods, the aforementioned result of Wilton has also been extended by Jaikin-Zapirain and the author \cite[Theorem A]{JaikinMorales24}, where we show that a finitely generated good group $G$ with $\cd(G)=2$ that is RFRS  (which stands for ``residually finite rationally solvable'' \cite{Ago13}) satisfies that its profinite completion is $\PD_2$  at some prime $p$ if only if $G$ is  a surface group itself.

In this article, we combine some well-known results on profinite completions of limit groups and three-manifold groups to study the recognition of  being a central extension of a surface (\cref{recog}). We also study the profinite invariance of coherence (\cref{coherence}) within the class of residually free groups. Another motivation of this work was extending a previous result of the author \cite[Theorem F]{morales24} to provide more  examples of groups that are recognised among finitely generated residually free groups by their finite quotients. This is ensured in \cref{rigid}, after having   establishing the following. 

\begin{thmA} \label{recog} Let $n\geq 0$ be an integer and let $M$ be a closed irreducible three-manifold. If $G$ is a finitely generated residually free group  with $\hat G\cong \hat{\pi_1M}\times \hat{\Z}^n$, then $G\cong H\times \Z^n,$ where $H$ is a central extension of a surface group $1\rar \Z\rar H\rar \pi_1 \Sg\rar 1$ and $\hat H \cong \hat{\pi_1M}.$
\end{thmA}

An ingredient needed in the proof of \cref{recog} is a criterion to read off when a central extension of groups splits from their induced extension of profinite completions (\cref{lem:exactgood}), which is inspired by the work of Wilton--Zalesskii  \cite{WilZal17} on distinguishing profinite completions of circle bundles over a surface. Some of the ideas involved in the proof of \cref{recog} also resonate in a recent work of Fruchter and the author \cite[Theorem E]{FruchterMorales24}, where we show that direct products of free and surface groups are profinitely rigid  within finitely presented residually free groups.

On a different note, Wilton  \cite{Wil08} classified which compact three-manifold groups with incompressible toral boundary have residually free fundamental group. In particular, he describes which groups $H$ satisfy  the conclusion of \cref{recog} as follows. 
\begin{prop}[Section 2.1, \cite{Wil08}]  \label{list} Let $H$ be a residually free group that fits in a central extension of the form $1\rar \Z\rar H\rar \pi_1 \Sg\rar 1$. Then one of the following happens.
\begin{itemize}
    \item $\Sg$ is orientable and  $H\cong \pi_1\Sg\times \Z$; or
    \item  $\Sg$ is non-orientable, with genus $g\geq 4$ and, for some $e\in \{0, 1\}$, \[H\cong \lan a_1, \dots, a_{g}, z\, |\, \prod_i a_i^2=z^e, [a_1, z]=\dots=[a_g, z]=1\ran. \]
\end{itemize}
\end{prop}

These constrains on the monodromy of $\Z$-by-surface groups as a consequence of the assumption on residual freeness become even more restrictive for free-by-cyclic groups (as reflected by \cref{freecyclic}). We now explain why \cref{recog} implies \cref{rigid}. On the one hand, Wilkes \cite[Theorem 6.1]{Wil17} classified which pairs of closed orientable Seifert-fibred spaces have the same profinite completion, which implies that each of the groups from \cref{list}  are profinitely rigid within the class of closed orientable Seifert-fibred spaces. Wilkes' classification was extended by Piwek \cite{Piw23} to the setting of central extensions of 2-orbifold groups with higher-rank centre. On the other hand, Wilton--Zalesskii \cite[Theorem B]{Wilton17} proved that being a Seifert-fibred space is a profinite invariant within closed aspherical three-manifolds, so the discussion of this paragraph together with \cref{recog} lead directly to the following.

\begin{corA} \label{rigid} Let $n\geq 0$ be an integer and let $\pi_1 M$ be one of the groups of \cref{list}. Suppose that  $G$ is a finitely generated residually free group with $\hat G\cong \hat {\pi_1 M}\times \hat{\Z }^n$. Then $G\cong \pi_1M\times \Z^n$.
\end{corA}

Regarding higher dimensional manifolds with residually free fundamental groups,  Fruchter and the author \cite[Corollary D]{FruchterMorales24} show that an aspherical closed topological manifold of dimension $n\geq 5$ has a finite-sheeted covering homeomorphic to a direct product of closed surfaces and (only when $n$ is odd) a circle. 

We move on to introduce \cref{freecyclic}. Two classes of two-dimensional groups that have received much attention in geometric group theory are one-relator groups and free-by-cyclic groups. On the one hand, we do not know which one-relator groups are residually free. Baumslag \cite{Bau67a} gave a method to construct many of these, although residual freeness still seems to be a very exceptional property. Under this additional assumption, Fruchter and the author confirm Melnikov's surface conjecture  \cite[Corollary D]{FruchterMorales24}. On the other hand, understanding which  free-by-cyclic groups are residually free is much simpler.
Bridson--Howie \cite[Corollary 3.4]{Bri07} showed that no group of the form (non-abelian finitely generated free)-by-$\Z$ is a limit group. We prove the following more general statement.

\begin{thmA} \label{freecyclic} Let $G$ be a finitely generated coherent  group that satisfies  $\b(G)=0$ and $\cd_p(\hat G)=2$. If $G$ is residually free, then there exists a free group $F$ such that $G\cong F\times \Z$.
\end{thmA}
A group $G$ of the form (finitely generated free)-by-$\Z$ has $\b(G)=0$ (see e.g. \cite[Theorem 1.39]{Luc02}) and it is coherent by Feighn--Handel \cite{Fei99}. Hence \cref{freecyclic} implies that such $G$ is residually free if only if it has the form $F\times \Z$.  

Our proof of \cref{recog} relies on the fact that a prime closed three-manifold $M$ has $\b(\pi_1 M)=0$ by  Lott--Lück \cite{Lot95}. When tackling the more general statement of \cref{recog} where $\pi_1 M$ is allowed to have a non-trivial prime decomposition, one encounters the hard problem of reading off free factors of a group from its profinite completion. This was solved in the context of three-manifold groups by Wilton--Zalesskii \cite{Wil19} and extended to a slightly bigger class by de Bessa--Porto--Zalesskii \cite{BessaPortoZalesskii2023}. However, to the best of the author's knowledge, this problem is open in greater generality. So \cref{recog} does not answer the following question. 
\begin{question} Suppose that $G$ is a residually free group and that $M$ is a three-manifold. Does $\hat G\cong \hat{\pi_1 M}$ imply that $G\cong \pi_1 M$?
\end{question}

Lastly, we note that coherence is a profinite invariant within residually free groups as part of a more general statement in \cref{coherence2}.
\begin{thmA}\label{coherence} Let $G$ and $H$ be two finitely generated residually free groups with $\hat G\cong \hat H$. Then $G$ is coherent if and only if $H$ is coherent. 
\end{thmA}
\subsection{Acknowledgements} 
The author is funded by the Oxford--Cocker Graduate Scholarship. This work has also received funding from the European Research Council (ERC) under the European Union's Horizon 2020 research and innovation programme (Grant agreement No. 850930). The author is grateful to Sam Fisher for a helpful conversation about \cref{lem:exactgood} and to an anonymous referee for various useful comments. Finally, the author would like to thank his supervisors Dawid Kielak and Richard Wade for their advice and support. 

\section{Residually free groups and the proof of \cref{coherence}}
We review several properties of residually free groups and several consequences which will play a role in the proofs of Theorems \ref{recog} and \ref{freecyclic}.  We start by listing a few elementary properties in the following lemma (the reader is referred to \cite{Bau62} for proofs).

\begin{lem} \label{props} Let $G$ be a finitely generated residually free group. 
\begin{enumerate}
    \item Every abelian subgroup of $G$ is finitely generated.
    \item The group $G/Z(G)$ is residually free. 
    \item If $K \leq Z(G)$ is a subgroup and $x\in G$  satisfies that $x K\in Z(G/K)$, then $x \in Z(G).$
    \item If $H\leq G$ is a finite-index subgroup, then $Z(H)=H\cap Z(G)$. 
\end{enumerate}
\end{lem}

The following relates the structure theory of residually free groups to that of limit groups. 

\begin{prop}[\cite{Bau99}, Corollary 19; and \cite{Sel01}, Claim 7.5]\label{subprod} Let $G$ be a finitely generated residually free group. Then $G$ is a subdirect product of finitely many limit groups.
\end{prop}

The conclusion of \cref{subprod} was pushed much further under additional finiteness assumptions by Bridson--Howie--Miller--Short as follows. 
\begin{prop}[Corollary 1.1, \cite{Bri09}] \label{vprod} Every $\FP_{\infty}(\Q)$ residually free group is virtually the direct product of finitely many limit groups. 
\end{prop}

Recall that a group $G$ is said to be $\FP_{\infty}(\Q)$ if the trivial $\Q[G]$-module $\Q$ admits a resolution by finitely generated projective modules (the definition of $\FP_{\infty}(\Z)$ that we refer to in \cref{ses} is completely analogous). We now collect a set of equivalences which we shall need later. They are undoubtedly well known to experts but we could not find a explicit formulation in the literature, so we  state it and prove it here for completeness. 

\begin{lem}\label{equiva} Let $G$ be a finitely generated residually free group. The following are equivalent: 
\begin{itemize}
    \item[$(a)$] There is no injection $F_2\times F_2\inc G$.
     \item[$(b)$]  $G$ is coherent. 
    \item[$(c)$] $G$ is LERF.
    \item[$(d)$]  $G$ is a central extension of a limit group, i.e. there exists a central abelian subgroup $A\leq G$ such that $G/A$ is a limit group.
    \item[$(e)$]  $G$ is virtually isomorphic to the direct product of a free abelian group and a limit group. 
\end{itemize}
\end{lem}
\begin{proof} The implications $(b)\implies (a)$ and $(c)\implies (a)$ are due to the fact that $F_2\times F_2$ is neither coherent nor LERF and that these properties pass to finitely generated subgroups. The implications $(d)\implies (b)$ and $(d)\implies (c)$ are due to the fact that limit groups are coherent and LERF \cite{Wil08b} and also the elementary fact that central extensions of a coherent (resp. LERF) group is also coherent (resp. LERF). Since $(e)$  clearly implies $(a)$, it remains to show  $(a)\implies (d)$ and $(d)\implies (e)$.

Let us assume $(a)$. By \cref{subprod}, $G$ is a subdirect product of finitely many limit groups $G\leq L_1'\times \cdots \times L_k'$. By the assumption of $(a)$, at most one such intersection is not an infinite cyclic group. From here it is not hard to see that we can recover $G$ from a limit group by extending finitely many times by a central cyclic subgroup. By the property (3) of \cref{props}, this process expresses $G$ as a central extension of a limit group, proving $(d)$. 

Now let us assume $(d)$. Since limit groups are $\FP_{\infty}(\Q)$, $G$ is $\FP_{\infty}(\Q)$ and, by \cref{vprod}, $G$ is virtually a direct product of limit groups $L_1\times\cdots \times L_n$. We also know that $G$ contains no copy of $F_2 \times F_2,$ so at most one of the factors $L_i$ is not a free abelian group. This proves the implication $(d) \implies (e)$. 
\end{proof}
We now prove \cref{coherence} from the introduction, which is integrated in the following result. 
\begin{thm}\label{coherence2} Let $G$ and $H$ be finitely generated residually free groups with $\hat G\cong \hat H$. Then $G$ satisfies any of the equivalent properties of \cref{equiva} if and only if $H$ does. 
\end{thm}
\begin{proof} Assume that $G$ contains a subgroup isomorphic to $F_2\times F_2$. Then, by \cite{Bri08}, $\hat G$ contains a closed subgroup isomorphic to $\hat F_2\times \hat F_2$. In particular, $\HH=(F_2)_{\hp}\times (F_2)_{\hp}\inc \hat H$, where $(F_2)_{\hp}$ denotes the pro-$p$ completion of $F_2$ (i.e. the free pro-$p$ group of rank two). We now proceed by contradiction and suppose that $H$ contains no $F_2\times F_2$. By \cref{equiva},  $H/Z(H)=L$ is a limit group. Since $\HH$ is centreless, and it is contained in $\hat H$, then $\HH$ is contained in $\hat L$. By a profinite Bass-Serre theoretical argument due to Zalesskii--Zapata \cite[Theorem C]{Zal20}, the pro-$p$ group $\HH$ is a free product of free abelian pro-$p$ groups. However, this  is a contradiction because $(F_2)_{\hp}\times (F_2)_{\hp}$ is not abelian and cannot possibly be a free product of non-trivial pro-$p$ groups. One way to show this last claim is by noting that, in such case,   \cite[Proposition 4.2.4]{Rib17} would imply that $\HH$ is projective. However, $\HH$ has $p$-cohomological dimension two.  This contradiction finishes the proof. 
\end{proof}

\section{Goodness and the proofs of Theorems \ref{recog} and \ref{freecyclic}} \label{ProofSection}
Before moving on to the proofs, we recall Serre's notion of goodness \cite{Ser97} and some of its fundamental properties. A group $G$ is said to be {\it good} if for all finite $G$-modules $A$ and all integers $n\geq 0$, the canonical maps $H^n(\hat G, A)\rar H^n(G, A)$, induced by the map $G\rar \hat G$, are isomorphisms. In particular, if $G$ is good, then $\cd_p(G)=\cd_p(\hat{G})$. We will repeatedly use the fact that if $1\rar K\rar G\rar Q\rar 1$ is a short exact sequence of groups, where $K$ is finitely generated and $Q$ is  good, then the induced sequence in profinite completions $1\rar \hat K\rar \hat G\rar \hat Q\rar 1$ is exact \cite[Ex 2(b), Section 2.5]{Ser97}.

The main sources of examples of good groups that we shall need either come from limit groups  \cite[Theorem 1.3]{Gru08} or from the ones that lie under the assumptions of the following standard lemma from \cite[Ex 2(c), Section 2.5]{Ser97} (see also \cite[Theorem 2.5]{Lor08} for a precise proof).

\begin{lem}\label{ses} Suppose that $1\rar K\rar G\rar Q\rar 1$ is a short exact sequence with $K, Q$ good and, in addition, $K$ of type $\FP_{\infty}(\Z)$. Then $G$ is good. 
\end{lem}

One additional lemma required in the proof of \cref{recog} describes a scenario where one can  profinitely recognise when a central extension split. 

\begin{lem}\label{lem:exactgood} Let $1\rar A\rar G\rar Q\rar 1$ be a short exact sequence of groups such that $A\leq Z(G)$, where $A$ is finitely generated and $Q$ is a finitely generated good group. If the induced sequence in profinite completions $1\rar\hat A\rar \hat G\rar \hat Q\rar 1$ splits, then so does the original one. 
\end{lem}
The previous lemma was discovered by Wilton--Zalesskii \cite[Lemma 8.3]{Wilton17}, who proved it when $Q$ is a surface group. They used this criterion to profinitely distinguish circle bundles over  a surface (and this was subsequently used in dimension 4 by Jiming--Zixi \cite[Proposition 25]{Ma22}). We give another proof that carries over to the current formulation.  \cref{lem:exactgood} is also employed   in the recognition of certain direct product decompositions among residually free groups in \cite[Theorem F]{FruchterMorales24}.

\begin{proof}[Proof of \cref{lem:exactgood}] We endow $A$ with the natural $Q$-module structure given by conjugation, which turns it into a trivial module. Given an integer $m\geq 1$, we denote by $A/m$ the quotient $A/mA$. The group $A/m$ is naturally a trivial $Q$-module and $A\rar A/m$ will denote the canonical reduction modulo $m$. We shall first prove  the following claim. 
\begin{claim} \label{claim:exact} Given a non-trivial element   $x\in H^2(Q; A)$, there exists an integer $m\geq 1$ such that $x$ does not belong to the kernel of the natural map $H^2(Q; A)\rar H^2(Q; A/m)$.     
\end{claim}
\begin{proof}[Proof of \cref{claim:exact}] Since $A$ is finitely generated, we can write $A\cong \Z^d\oplus U$ for some integer $d\geq 0$ and a finite abelian subgroup $U$. Since we have the natural splitting $H^2(Q; A)\cong H^2(Q; \Z)^d\oplus H^2(Q; U)$, it directly follows that we only need to prove the claim when $A\cong \Z$. By the universal coefficient theorem for cohomology \cite[Theorem 3.2]{Hat00}, we have the following commutative diagram of natural maps  
\begin{equation*}
    \begin{tikzcd}[cramped, sep=small]
       0\ar[r] &  \Ext_{\Z}^1(H_1(Q; \Z), \Z)\ar[d, "\pi_1^m"] \ar[r, "\alpha"]& H^2(Q; \Z) \ar[d, "\pi_2^m"]\ar[r, "\beta"]& \Hom_{\Z}(H_2(Q; \Z), \Z)\ar[d, "\pi_3^m"]\ar[r]&   0\\ 
       0\ar[r] &  \Ext_{\Z}^1(H_1(Q; \Z), \Z/m) \ar[r]& H^2(Q; \Z/m) \ar[r]& \Hom_{\Z}(H_2(Q; \Z), \Z/m)\ar[r]&   0. \end{tikzcd}
\end{equation*}
An important, preliminary observation is that the torsion of $H^2(Q;  \Z)$ comes from torsion in the abelianisation of $Q$. To see this, we note that the group $\Hom_{\Z}(H_2(Q; \Z), \Z)$ is  torsion-free.  On the other hand, if $H_1(Q;  \Z)\cong \Z^r\oplus T$, where $T$ is a finite abelian group and $r\geq 0$ is an integer; then $\Ext_{\Z}^1(H_1(Q; \Z), \Z)\cong T$ and $\Ext_{\Z}^1(H_1(Q; \Z), \Z/m)\cong T/m$ (these computations are explained right after the statement of \cite[Theorem 3.2]{Hat00}). Furthermore, the map $\pi_1^m$ from the diagram above corresponds, under the previous isomorphisms, to the natural  map  $T\rar T/m$ of reduction modulo $m$. We are now in the right position to prove \cref{claim:exact}. If $x\in H^2(Q; \Z)$ is not torsion, then $\beta(x)$ is a non-trivial element of $\Hom_{\Z}(H_2(Q; \Z), \Z)$. So it is clear that $\pi_3^m(\beta(x))\neq 0$. In particular, $\pi_2^m(x)\neq 0$ in this case. Finally, it remains to reach the same conclusion when $x$ is torsion. In this case, $x\in \ker \beta=\im \alpha$. Furthermore, by our discussion above, if we take $m=|T|$, then $T/m\cong T$ and $\pi_1^m(\alpha^{-1}(x))\neq 0$ (since $\pi_1^m$ is an isomorphism). The proof is complete.   \renewcommand\qedsymbol{$\diamond$}
\end{proof} Once \cref{claim:exact} is proved, we can argue as in \cite[Lemma 8.3]{Wilton17}. Let $B$ denote a finite $Q$-module, which can be naturally be viewed as a $\hat{Q}$-module. It is well-known that $H^2(Q; B)$ is naturally identified with the equivalence classes of extensions of $B$ by $Q$ (see, for instance, \cite[Chapter IV, Section 3]{Bro82}). An analogous principle carries over to the setting of profinite groups regarding the cohomology group $H^2(\hat Q; B)$ (see \cite[Section 6.8]{Rib00}). In this sense, the natural map $H^2(\hat Q; B)\rar H^2(Q; B)$ induced by restriction can be realised in terms of pulling back   extensions of $B$ by $\hat Q$ (say, of the form $1\rar B\rar \Gamma\rar \hat Q\rar 1$ for some profinite group $\Gamma$) to   extensions of $B$ by $Q$ (of the form $1\rar B\rar P\rar Q\rar 1$). This pull-back construction is explained in detail in Wilkes' thesis \cite[Proposition 2.2.4]{Wilkes2018} and an important feature of it is that, in our previous notation, there is a relation between the middle terms $\Gamma$ and $P$ of the extensions given by a natural isomorphism $\Gamma\cong \hat P$.

Now we are ready to complete the proof of \cref{lem:exactgood}. Suppose that $1 \rar A\rar G\rar Q\rar 1$ does not split. Then it corresponds to some non-zero element $x\in H^2(Q; A)$. By \cref{claim:exact}, there exists an integer $m\geq 1$ such that the natural sequence $1\rar A/mA\rar G/mA\rar Q\rar 1$ does not split. By the goodness of $Q$, the map $H^2(\hat Q; A/mA)\rar H^2(Q; A/mA)$ is an isomorphism and the sequence $1\rar A/mA\rar G/mA\rar Q\rar 1$ will be obtained from the pull-back of the (necessarily non-split) sequence $1\rar A/mA\rar \hat G/m\hat A\rar \hat Q\rar 1$. In particular, the sequence $1\rar \hat A\rar \hat G \rar \hat Q\rar 1$ does not split, as we wanted. This finishes the proof of \cref{lem:exactgood}.
\end{proof} 
The   ideas that appear in the proof of \cref{lem:exactgood} are  similar to those from \cite[Proposition 6.1]{Gru08} on ensuring that finite extensions of residually finite good groups are residually finite.

We come back to discussing the proof of \cref{recog} from the introduction, but we shall first recall its statement. 
\begin{thm}[Theorem A] \label{recog2} Let $n\geq 0$ be an integer and let $M$ be a closed irreducible three-manifold. If $G$ is a finitely generated residually free group $G$ with $\hat G\cong \hat{\pi_1M}\times \hat{\Z}^n$, then $G\cong H\times \Z^n,$ where $H$ is a central extension of a surface group $1\rar \Z\rar H\rar \pi_1 \Sg\rar 1.$
\end{thm}

The proof of \cref{recog2} is split into several claims. If $\pi_1 M$ is abelian, then \cref{recog2} is trivial. So we assume from now on that $\pi_1 M$ and $G$ are not abelian.

\begin{claim}\label{F2} The group $G$ contains no subgroup isomorphic to $F_2\times F_2$. 
\end{claim}
\begin{proof} We proceed by contradiction. Suppose that there is an injection $F_2\times F_2\inc G$. The group $F_2\times F_2$ is $\FP_{\infty}(\Q)$ (and so are all its finite-index subgroups), and all these groups are proven to be separable in $G$ by Bridson--Wilton \cite{Bri08}. So $G$ induces the full profinite topology on $F_2\times F_2$ and the induced map $\hat F_2\times \hat F_2\inc \hat G$ is injective. This implies that the pro-$p$ group $\HH=(F_2)_{\hp}\times (F_2)_{\hp}$ is a closed subgroup of $\hat G$. Since $\HH$ is centreless, we get in particular that $\HH\inc \hat{\pi_1 M}$ for all primes $p$. The pro-$p$ subgroups of $\hat{\pi_1 M}$ are classified by Wilton--Zalesskii \cite[Theorem 1.3]{WilZal17} and these are, for $p>3$, coproducts of groups of the following list: $\Z/p$, $\Z_p$, $\Z_p^2$, pro-$p$ completions of $(\Z\times \Z)\rtimes \Z$ and the pro-$p$ completion of a residually-$p$ fundamental group of a non-compact Seifert-fibred manifold with hyperbolic base orbifold. This leads to a contradiction because $\HH$ does not belong to this list and, furthermore, it cannot be the coproduct of two non-trivial pro-$p$ groups (arguing exactly as in the proof of \cref{coherence2}).
\end{proof}

\begin{claim} \label{extension} The group $G$ is a central extension $1\rar \Z^{n+1}\rar G\rar L\rar 1$ where $L$ is a non-abelian limit group.
\end{claim}
\begin{proof} From \cref{F2} and \cref{equiva}, we deduce that the group $G$ is a central extension $1\rar \Z^{m}\rar G\rar L\rar 1$ for some $m\geq 0$ and some limit group $L$. By property (3) of \cref{props}, $L$ must be non-abelian (we assumed $G$ to be non-abelian). We will prove that $m=n+1$. The group $L$ is good  and hence there is an induced exact sequence on profinite completions $1\rar \hat{\Z}^m\rar \hat G\rar \hat L\rar 1$. By  \cite[Corollary 4.4]{Zal20}, $\hat L$ is centreless and hence $Z(\hat G)\cong \hat{\Z}^m$. If the centre of $\hat{\pi_1 M} $ was trivial, then $\hat{\pi_1 M}\cong \hat G/Z(\hat G)\cong \hat L$, which is absurd because from Lück's approximation \cite{Luc94} we would derive  that $0=\b(\pi_1 M)=\b(L)$, contradicting the fact that $L$ is a non-abelian limit group. This proves that the centre of $\hat{\pi_1 M} $ is non-trivial and hence $m\geq n+1$. If $m\geq n+2$, then \cite[Theorem 3.7.4]{Neu08} would yield that $\hat{L}$ is $\PD_{n+3-m}$ at $p$ (i.e. it is $\PD_0$ or $\PD_1$), implying that $L$ is cyclic, which is a contradiction. 
\end{proof}

\begin{claim} \label{extension2} The group $G$ is isomorphic to $\Z^n\times G_1$, where $\hat{G_1}\cong \hat{\pi_1M}$ and $G_1$ fits in a central extension of the form $1\rar \Z\rar G_1\rar \pi_1\Sg\rar 1$ for some closed  surface $\Sg$. 
\end{claim}
\begin{proof} The proof of \cref{extension} shows that $\ove{Z(G)}=Z(\hat{G})$, so $Z(G)$ is dense in $Z(\hat G).$ By assumption, we have a decomposition $\hat G\cong \hat{\pi_1M}\times \hat{\Z}^n$. We denote by $\KK$ the closed subgroup corresponding to the direct factor $\hat{\Z}^n$ of $\hat G$. There exists $K\leq Z(G)$ isomorphic to $\Z^n$ such that $\ove K=\KK$. By \cref{extension}, there is a short exact sequence $1\rar \Z\rar G/K\rar L\rar 1$, for a limit group $L$, so $G/K$ is good by \cref{ses}. Since the  short exact sequence in profinite completions $1\rar \KK \rar \hat G\rar \hat{ G/K}\rar 1$ splits, it follows from \cite[Lemma 8.3]{Wilton17}  that the short exact sequence $1\rar K\rar G\rar G/K\rar 1$ splits and hence $G\cong K\times G/K$, where $\hat{G/K}\cong  \hat{\pi_1M}$.   In particular, $\hat{G/K}$ is $\PD_3$ at all primes $p$ by \cite[Theorem 4.1]{Koc08}. From the short exact sequence $1\rar \hat \Z\rar \hat{G/K}\rar \hat L\rar 1$ and \cite[Theorem 3.7.4]{Neu08}, we also derive that $\hat L$ is $\PD_2$ at all primes $p$. By an argument of Wilton \cite{Wil21}, this implies that $L$ is a surface group. Regarding the latter implication, we note that Wilton's assumptions are that $\hat L$ is the profinite completion of a surface group, but his arguments only use the {\it a priori} weaker fact that $\hat L$ is $\PD_2$ at some prime $p$ (see \cite[Proposition 8.10]{morales24} for a more detailed account of this argument). 
\end{proof}

This completes the proof of \cref{recog2}. Lastly, we restate and prove   \cref{freecyclic} from the introduction.

\begin{thm}[\cref{freecyclic}] Let $G$ be a finitely generated coherent  group that satisfies $\b(G)=0$ and $\cd_p(\hat G)=2$. If $G$ is residually free, then there exists a free group $F$ such that $G\cong F\times \Z$.
\end{thm}
\begin{proof} By \cref{equiva}, there is a central extension $1\rar \Z^m\rar G\rar L\rar 1$, where $L$ is a limit group. The group $L$ is good and so is $G$ by \cref{ses}. By the (virtual) additivity of $\cd_p$ established by Weigel--Zalesskii \cite[Theorem 1.1]{Wei04}, $m+\cd_p(\hat L)=\cd_p(\hat G)=2$. If $G$ is abelian, it is clear that $G\cong \Z^2$. If $G$ is non-abelian, then $L\neq 1$ and $\cd_p(L)=\cd_p(\hat L)\geq 1$, so $m\leq 1$. In addition, the fact that $G$ is non-abelian and $\b(G)=0$ implies that $G$ cannot be a limit group, so $m\geq 1$. Passing to profinite completions, we get an exact sequence  $1\rar \hat{\Z}\rar \hat{G}\rar \hat{L}\rar 1$. This leads to  $\cd_p(\hat L)=1$, which implies that $L$ is free \cite[Theorem 7.1]{Wil18}. So the central extension $1\rar \Z\rar G\rar L\rar 1$ splits and we are done. 
\end{proof}

\textsc{Mathematical Institute, University of Oxford, Radcliffe Observatory, Andrew Wiles Building, Woodstock Rd, Oxford OX2 6GG} \\
\textit{E-mail address:} \href{mailto:morales@maths.ox.ac.uk}{\texttt{morales@maths.ox.ac.uk}}

\bibliography{biblio.bib}

\begin{bibdiv}
\begin{biblist}

\bib{Ago13}{article}{
      author={Agol, Ian},
       title={The virtual haken conjecture (with an appendix by ian agol,
  daniel groves and jason manning)},
        date={2013},
     journal={Doc. Math.},
      volume={18},
       pages={1045\ndash 1087},
}

\bib{Bau62}{article}{
      author={Baumslag, G},
       title={On generalised free products},
        date={1962},
     journal={Mathematische Zeitschrlft},
      volume={78},
       pages={423\ndash –438},
}

\bib{Bau67a}{article}{
      author={Baumslag, Benjamin},
       title={Residually free groups},
        date={1967},
     journal={Proceedings of the London Mathematical Society},
      volume={s3-17},
      number={3},
       pages={402\ndash 418},
}

\bib{Bri07}{article}{
      author={Bridson, Martin},
      author={Howie, James},
       title={Normalisers in limit groups},
        date={2007},
     journal={Mathematische Annalen},
      volume={337},
       pages={385\ndash 394},
}

\bib{Bri09}{article}{
      author={Bridson, Martin~R.},
      author={Howie, James},
      author={Miller, Charles~F.},
      author={Short, Hamish},
       title={Subgroups of direct products of limit groups},
        date={2009},
     journal={Annals of Mathematics},
      volume={170},
      number={3},
       pages={1447\ndash 1467},
}

\bib{Bau99}{article}{
      author={Baumslag, Gilbert},
      author={Myasnikov, Alexei},
      author={Remeslennikov, Vladimir},
       title={Algebraic geometry over groups},
        date={1999},
     journal={Journal of algebra},
      volume={219},
       pages={16\ndash 79},
}

\bib{Bro82}{book}{
      author={Brown, Kenneth~S.},
       title={Cohomology of groups},
      series={Graduate Texts in Mathematics},
   publisher={Springer, New York, NY},
        date={1982},
}

\bib{Bri08}{article}{
      author={Bridson, M.~R.},
      author={Wilton, H.},
       title={Subgroup separability in residually free groups},
        date={2008},
     journal={Mathematische Zeitschrift},
      volume={260},
       pages={25\ndash 30},
}

\bib{BessaPortoZalesskii2023}{article}{
      author={de~Bessa, Vagner~R.},
      author={Porto, Anderson L.~P.},
      author={Zalesskii, Pavel~A.},
       title={The profinite completion of accessible groups},
        date={2023},
        ISSN={0026-9255,1436-5081},
     journal={Monatsh. Math.},
      volume={202},
      number={2},
       pages={217\ndash 227},
         url={https://doi.org/10.1007/s00605-022-01789-9},
      review={\MR{4641682}},
}

\bib{Fei99}{article}{
      author={Feighn, Mark},
      author={Handel, Michael},
       title={Mapping tori of free group automorphisms are coherent},
        date={1999},
     journal={Annals of Mathematics. Second Series},
      volume={149},
      number={3},
       pages={1061\ndash 1077},
}

\bib{FruchterMorales24}{misc}{
      author={Fruchter, Jonathan},
      author={Morales, Ismael},
       title={Virtual homology of limit groups and profinite rigidity of direct
  products},
        date={2024},
        note={arxiv:2209.14925},
}

\bib{Gru08}{article}{
      author={Grunewald, F.},
      author={Jaikin-Zapirain, A.},
      author={Zalesskii, P.~A.},
       title={Cohomological goodness and the profinite completion of bianchi
  groups},
        date={2008},
     journal={Duke Mathematical Journal},
      volume={144},
      number={1},
       pages={53\ndash –72},
}

\bib{Hat00}{book}{
      author={Hatcher, Allen},
       title={{Algebraic topology}},
   publisher={Cambridge University Press},
        date={2000},
}

\bib{JaikinMorales24}{misc}{
      author={Jaikin-Zapirain, Andrei},
      author={Morales, Ismael},
       title={Prosolvable rigidity of surface groups},
        date={2024},
        note={arxiv:2312.12293},
}

\bib{Koc08}{article}{
      author={Kochloukova, Dessislava~H.},
      author={Zalesskii, Pavel~A.},
       title={Profinite and pro-p completions of poincaré duality groups of
  dimension 3},
        date={2008},
     journal={Transactions of the American Mathematical Society},
      volume={360},
      number={4},
       pages={1927\ndash 1949},
}

\bib{Luc02}{book}{
      author={Lück, Wolfgang},
       title={{$L^2$}-{I}nvariants: {T}heory and {A}pplications to {G}eometry
  and {K}-theory},
      series={Ergebnisse der Mathematik und ihrer Grenzgebiete. 3. Folge / A
  Series of Modern Surveys in Mathematics},
   publisher={Springer Berlin, Heidelberg},
        date={2002},
}

\bib{Luc94}{article}{
      author={Lück, W.},
       title={Approximating {$L^2$}-invariants by their finite-dimensional
  analogues},
        date={1994},
     journal={Geometric and Functional Analysis},
       pages={455\ndash –481},
}

\bib{Lot95}{article}{
      author={Lott, John},
      author={L\"{u}ck, Wolfgang},
       title={{$L^2$}-topological invariants of {$3$}-manifolds},
        date={1995},
     journal={Inventiones Mathematicae},
      volume={120},
      number={1},
       pages={15\ndash 60},
}

\bib{Lor08}{article}{
      author={Lorensen, Karl},
       title={Groups with the same cohomology as their profinite completions},
        date={2008},
     journal={Journal of Algebra},
      volume={320},
      number={4},
       pages={1704\ndash 1722},
}

\bib{morales24}{article}{
      author={Morales, Ismael},
       title={On the profinite rigidity of free and surface groups},
        date={2024},
        ISSN={0025-5831,1432-1807},
     journal={Math. Ann.},
      volume={390},
      number={1},
       pages={1507\ndash 1540},
         url={https://doi.org/10.1007/s00208-023-02785-6},
      review={\MR{4800944}},
}

\bib{Ma22}{article}{
      author={Ma, Jiming},
      author={Wang, Zixi},
       title={Distinguishing 4-dimensional geometries via profinite
  completions},
        date={2022},
     journal={Geometriae Dedicata},
      volume={216},
      number={52},
}

\bib{Neu08}{book}{
      author={Neukirch, J{\"u}rgen},
      author={Schmidt, Alexander},
      author={Wingberg, Kay},
       title={Cohomology of number fields},
      series={Grundlehren der mathematischen Wissenschaften},
   publisher={Springer Berlin, Heidelberg},
        date={2008},
}

\bib{Piw23}{misc}{
      author={Piwek, Pawel},
       title={Profinite rigidity properties of central extensions of 2-orbifold
  groups},
        date={2023},
        note={arxiv:2304.01105},
}

\bib{Rib17}{book}{
      author={Ribes, Luis},
       title={Profinite graphs and groups},
      series={Ergebnisse der Mathematik und ihrer Grenzgebiete. 3. Folge. A
  Series of Modern Surveys in Mathematics [Results in Mathematics and Related
  Areas. 3rd Series. A Series of Modern Surveys in Mathematics]},
   publisher={Springer, Cham},
        date={2017},
      volume={66},
}

\bib{Rib00}{book}{
      author={Ribes, Luis},
      author={Zalesskii, Pavel},
       title={Profinite groups},
   publisher={Springer Berlin Heidelberg},
        date={2000},
}

\bib{Sel01}{article}{
      author={Sela, Zlil},
       title={Diophantine geometry over groups {I}: Makanin-razborov diagrams},
        date={2001},
     journal={Publications Math\'ematiques de l'IH\'eS},
      volume={93},
       pages={31\ndash 105},
}

\bib{Ser97}{book}{
      author={Serre, Jean-Pierre},
       title={Galois cohomology},
      series={Springer Monographs in Mathematics},
   publisher={Springer, Berlin, Heidelberg},
        date={1997},
}

\bib{Wil08b}{article}{
      author={Wilton, Henry},
       title={Hall's theorem for limit groups},
        date={2008},
     journal={Geometric and Functional Analysis},
      volume={18},
       pages={271\ndash –303},
}

\bib{Wil08}{article}{
      author={Wilton, Henry},
       title={{Residually free $3$–manifolds}},
        date={2008},
     journal={Algebraic {\&} Geometric Topology},
      volume={8},
      number={4},
       pages={2031 \ndash  2047},
}

\bib{Wil17}{article}{
      author={Wilkes, Gareth},
       title={Profinite rigidity for seifert fibre spaces},
        date={2017},
     journal={Geometriae Dedicata},
      volume={188},
       pages={141\ndash 163},
}

\bib{Wilkes2018}{article}{
      author={Wilkes, Gareth},
       title={Profinite properties of {3-Manifold} groups},
        date={2018},
        note={Dphil thesis},
}

\bib{Wil18}{article}{
      author={Wilton, Henry},
       title={Essential surfaces in graph pairs},
        date={2018},
     journal={Journal of the American Mathematical Society},
      volume={31},
       pages={893–919},
}

\bib{Wil21}{article}{
      author={Wilton, Henry},
       title={On the profinite rigidity of surface groups and surface words},
        date={2021},
     journal={Comptes Rendus. Math\'ematique},
      volume={359},
      number={2},
       pages={119\ndash 122},
}

\bib{Wei04}{article}{
      author={Weigel, Thomas~S.},
      author={Zalesskii, Pavel},
       title={Profinite groups of finite cohomological dimension},
        date={2004},
     journal={Comptes Rendus Mathematique},
      volume={338},
       pages={353\ndash 358},
}

\bib{Wilton17}{article}{
      author={Wilton, Henry},
      author={Zalesskii, Pavel},
       title={Distinguishing geometries using finite quotients},
        date={2017},
     journal={Geometry \& Topology},
      volume={21},
       pages={345\ndash 384},
}

\bib{WilZal17}{article}{
      author={Wilton, Henry},
      author={Zalesskii, Pavel},
       title={Pro-{$p$} subgroups of profinite completions of 3-manifold
  groups},
        date={2017},
     journal={Journal of the London Mathematical Society. Second Series},
      volume={96},
      number={2},
       pages={293\ndash 308},
}

\bib{Wil19}{article}{
      author={Wilton, Henry},
      author={Zalesskii, Pavel},
       title={Profinite detection of 3-manifold decompositions},
        date={2019},
     journal={Compositio Mathematica},
      volume={155},
      number={2},
       pages={246–259},
}

\bib{Zal20}{article}{
      author={Zalesskii, Pavel},
      author={Zapata, Theo},
       title={Profinite extensions of centralizers and the profinite completion
  of limit groups},
        date={2020},
     journal={Revista Matem\'{a}tica Iberoamericana},
      volume={36},
      number={1},
       pages={61\ndash 78},
}

\end{biblist}
\end{bibdiv}

\end{document}